\documentclass{amsproc}
\usepackage{graphicx}
\usepackage{amssymb}
\usepackage{epstopdf}
\usepackage{amsmath,amsthm,amscd,amssymb}
\usepackage{latexsym}
\usepackage[colorlinks,citecolor=red,pagebackref,hypertexnames=false]{hyperref}
\usepackage{geometry}                % See geometry.pdf to learn the layout options. There are lots.
\numberwithin{equation}{section}

\theoremstyle{plain}
\newtheorem{theorem}{Theorem}[section]
\newtheorem{lemma}[theorem]{Lemma}
\newtheorem{corollary}[theorem]{Corollary}

\theoremstyle{definition}
\newtheorem{definition}[theorem]{Definition}

\newtheorem{problem}[theorem]{Problem}
\newtheorem{example}[theorem]{Example}

\newtheorem{case[theorem]}{Case}

\def\vx{\vec{x}}

\def\vu{\vec{u}}
\def\vxi{\vec\xi}
\def\veta{\vec\eta}
\def\vzeta{\vec\zeta}
\def\R{\mathbb R}

\theoremstyle{remark}
\newtheorem{remark}[theorem]{Remark}

\numberwithin{equation}{section}

\begin{document}

\title[Restricted convolutions and multilinear operators]{\parbox{14cm}{\centering{Restricted convolution inequalities, multilinear operators \\ and applications}}}

\author{Dan-Andrei Geba, Allan Greenleaf, Alex Iosevich,  \\ Eyvindur Palsson and Eric Sawyer} 

\address{Department of Mathematics \\ University of Rochester\\ Rochester, NY 14627}
\email{dangeba@math.rochester.edu}
%\address{Department of Mathematics \\ University of Rochester\\ Rochester, NY 14627}
\email{allan@math.rochester.edu}
%\address{Department of Mathematics \\ University of Rochester\\ Rochester, NY 14627}
\email{iosevich@math.rochester.edu}
%\address{Department of Mathematics \\ University of Rochester\\ Rochester, NY 14627}
\email{palsson@math.rochester.edu}
\address{Department of Mathematics and Statistics, McMaster University,
Hamilton, ON L8S 4K1}
\email{sawyer@mcmaster.ca}

\thanks{The first three authors were supported by NSF grants DMS-0747656, DMS-0853892 and DMS-1045404, resp.  }

\begin{abstract} For $ 1\le k <n$, we prove that for functions $F,G$ on $ {\Bbb R}^{n}$, any $k$-dimensional affine subspace  $H \subset {\Bbb R}^{n}$, 
and $p,q,r \ge 2$ with $\frac{1}{p}+\frac{1}{q}+\frac{1}{r}=1$, one has the estimate

$$ {||(F*G)|_H||}_{L^{r}(H)} \leq {||F||}_{\Lambda^H_{2, p}({\Bbb R}^{n})} \cdot {||G||}_{\Lambda^H_{2, q}({\Bbb R}^{n})},$$
where the mixed norms on the right are defined by

$$ {||F||}_{\Lambda^H_{2,p}({\Bbb R}^{n})}={\left( \int_{H^*} {\left( \int {|\widehat{F}|}^2 
dH_{\xi}^{\perp} \right)}^{\frac{p}{2}} d\xi \right)}^{\frac{1}{p}},$$ 
with $dH_{\xi}^{\perp}$ the $(n-k)$-dimensional  Lebesgue measure on the affine subspace $H_{\xi}^{\perp}:=\xi + H^\perp$. 
Dually, one obtains  restriction theorems for the Fourier transform for affine subspaces.
Applied to  $F(x^{1},\dots,x^{m})=\prod_{j=1}^m f_j(x^{j})$ on $\R^{md}$, the diagonal $H_0=\{(x,\dots,x): x \in {\Bbb R}^d\}$ and suitable kernels $G$, this implies  new results for multilinear convolution operators, including $L^p$-improving bounds for  measures, an $m$-linear variant of Stein's spherical maximal theorem, estimates for $m$-linear oscillatory integral operators,  certain Sobolev trace inequalities, and bilinear estimates for solutions to the wave equation.

\end{abstract} 

\maketitle

\tableofcontents

\section{Introduction} 

\vskip.125in 

Convolution inequalities play a central role in harmonic analysis and related areas. Perhaps the most classical of these is Young's inequality:  if $F,G \in {\mathcal S}({\Bbb R}^n)$, then, for  $p,q,r \ge 1$,
\begin{equation} \label{young} 
{||F*G||}_{L^r({\Bbb R}^n)} \leq {||F||}_{L^p({\Bbb R}^n)} {||G||}_{L^q({\Bbb R}^n)},\quad \frac{1}{r}=\frac{1}{p}+\frac{1}{q}-1. \end{equation} 
Many of the applications of this inequality take the following form. (See, e.g., \cite{St93,Be75} and the references  there.) 
Let $F(x)=f(x)$, a function, and $g(x)=K(x)$, a suitable kernel or  measure. 
Then  (\ref{young})  implies Lebesgue spaces mapping properties of the linear convolution operator, 

$$Tf(x)=\int f(x-y) K(y) dy.$$

In this paper we prove an analogue of Young's inequality,  playing a similar role for \emph{multilinear} convolution operators as Young's inequality plays for linear ones.  This will follow from a result of independent interest, concerning restrictions of convolutions to affine subspaces; to state that, we need to introduce some notation.

For $1\le k\le n$, let $H\subset\Bbb R^{n}$ be a $k$-dimensional  linear subspace, $i_H:H\hookrightarrow \Bbb R^{n}$ the inclusion map, and $i_H^*:\Bbb (R^{n})^*\longrightarrow H^*$ the dual restriction map. $(\Bbb R^{n})^*$ foliates into a union of $(n-k)$-dimensional affine  subspaces, parallel to $(i_H^*)^{-1}(0)=H^\perp$,
\begin{equation}\label{foliate}
(\Bbb R^{n})^*=\bigcup_{\xi\in H^*} (i_H^*)^{-1}(\xi).
\end{equation}
Via the Euclidean metric, one identifies $(i_H^*)^{-1}(\xi)$ with $H_{\xi}^{\perp}:= \xi + H^\perp$.
The Lebesgue measure $d\zeta$ on $(\Bbb R^{n})^*$ thus decomposes as $d\zeta= d H_{\xi}^{\perp}(\zeta) \, d\xi$, where $d H_{\xi}^{\perp}$ and $d\xi$, resp., are the Lebesgue measures of dimensions $n-k$   and $k$  on $ H_{\xi}^{\perp}$ and $H^*$, resp.
Relative to this decomposition of the Fourier variables, we define a scale of mixed-norm spaces:

\begin{definition}\label{def space}
For $1\le r,p\le\infty$ and $F\in\mathcal S(\Bbb R^{n})$, let

\begin{equation} \label{Lambda}
{||F||}_{\Lambda^H_{r,p}({\Bbb R}^{n})}={\left( \int_{H^*} {\left[ \int_{H_{\xi}^{\perp}} {\left|\widehat{F}(\zeta)\right|}^r dH_{\xi}^{\perp}(\zeta) \right]}^{\frac{p}{r}} d\xi \right)}^{\frac{1}{p}},
\end{equation} 
and $\Lambda^H_{r,p}:=\Lambda^H_{r,p}({\Bbb R}^{n})$ be the completion of $\mathcal S({\Bbb R}^{n})$ with respect to this norm.

\end{definition}

\begin{remark} If $k=n$ or $r=p$, $||F||_{\Lambda^H_{r,p}}=||\widehat{F}||_{L^p}$, and if $r=p=2$, this equals $||F||_{L^2}$.
Furthermore, the norm is translation-invariant, so it is natural, if $\tilde{H}\subset\R^{n}$ is an \emph{affine} $k$-plane, to define $||\cdot||_{\Lambda^{\tilde{H}}_{r,p}}$ by using $H$ and $H_\xi^\perp$ on the right hand side of (\ref{Lambda}), where $H$ is the translate of $\tilde{H}$ passing through 0.
\end{remark}

We  now state our first result, which concerns restrictions of convolutions to affine subspaces.

\begin{theorem} \label{bilinearnewwithattitude} Let $ 1\le k<n$, $H\subset {\Bbb R}^{n}$ a $k $-dimensional affine   subspace,
and $F,G \in {\mathcal S}({\Bbb R}^{n})$. 
Then, for  $\frac{1}{p}+\frac{1}{q}+\frac{1}{r}=1$, $p,q,r \ge 2$,

\begin{equation} \label{newyoung} {||(F*G)|_H||}_{L^{r}(H)} \leq {||F||}_{\Lambda^{H}_{2, p}({\Bbb R}^{n})} \cdot {||G||}_{\Lambda^{H}_{2, q}({\Bbb R}^{n})}, \end{equation} 

\vskip.125in 

\noindent Moreover, the best constant in  (\ref{newyoung}), considering $F\rightarrow F*G|_H$  as a linear operator acting on $F$, is
\begin{equation} \label{bestconstant}
 {||F\longrightarrow F*G||}_{L^2(\Bbb R^{n})\longrightarrow L^2(H)} =
  C:= {||G||}_{\Lambda^{H}_{2,\infty}({\Bbb R}^{n})}. \end{equation} 

\end{theorem}

\vskip.125in

One motivation for considering  restricted convolution inequalities such as  (\ref{newyoung}) is the study of \emph{multilinear} convolution operators. For $d\ge 1,\, m\ge 2$, denote elements of $\R^{md}$ by $\vx:=(x^{1},\dots,x^{m})$ and elements of $(\R^{md})^*$ by $\vxi:=(\xi^{1},\dots,\xi^{m})$. 
\footnote{The notation $\vx,\, \vxi$, etc., will occasionally also be used for  $l$-tuples of vectors in $\R^d$  for $l<m$.}
For a kernel  $ K(\vx )\in \mathcal D'(\R^{md})$, define an $m$-linear convolution operator $T: \big(\mathcal D(\R^d)\big)^m\to \mathcal  E(\R^d)$,

\begin{equation}\label{multi T}
T(f_1,\dots,f_m)(x)=\int_{\R^{md}} f_1(x-u^1)\cdots f_m( x-u^m)\,  K(\vu)\, du^1\cdots du^m,\quad x\in\R^d.
\end{equation}
 
\noindent Then, since 

\begin{eqnarray*}\nonumber
T(f_1,\dots,f_m)(x)&=& \left(\left(f_1\otimes\cdots\otimes f_m\right)*K\right)(x,\dots,x),
\end{eqnarray*}
 multilinear bounds for $T$ can be obtained from Thm. \ref{bilinearnewwithattitude}, with $H$ being the \emph{diagonal} of $\R^{md}$, $H_0:=\{(x,\dots,x): x\in\R^d\}$. An immediate consequence is:
 
 \vskip.125in

 \begin{corollary} \label{bilinearnew} Let  $K\in\mathcal E'(\R^{md})$,  
 $f_1,\dots,f_m \in {\mathcal S}({\Bbb R}^d)$,
and $2 \leq r \leq \infty$.
Then,  the multilinear convolution operator defined by (\ref{multi T}) satisfies 

\begin{equation} \label{L2new} 
{||T(f_1,\dots,f_m)||}_{L^r(\R^d)} \lesssim ||K||_{\Lambda^{H_0}_{2,\frac{2r}{r-2}}}
\cdot\prod_{j=1}^m {||f_j||}_{L^2(\R^d)}.
\end{equation} 
\end{corollary} 
\vskip.125in 

\begin{remark} 
In Sec. \nolinebreak\ref{bilinearapplications}, we  show how Cor. \ref{bilinearnew} and some additional machinery can applied to obtain   $L^p$-improving results for multilinear convolution with measures, analogous to those in the linear setting. 
\end{remark}

 Note that $H_0^*$ can be identified with $\big\{\big(\frac{\xi}m,\dots,\frac{\xi}m\big): \xi\in (\R^d)^*\big\}$, while on  the \emph{anti-diagonal}, $H_0^\perp=\{\veta=(\eta^1,\dots,\eta^m):\, \sum\eta^j=0\}$,  we can solve for $\eta^m$ in terms of $(\eta^1,\dots,\eta^{m-1})$. This give rise to a parametrization making concrete the decomposition (\ref{foliate}),
 
 \begin{equation}\label{rhoH}
 \rho_{H_0}:\R^d_\xi\times\R^{(m-1)d}_{\veta}\to (\R^{md})^*,\quad \rho_H(\xi,\veta)=\big(\frac{\xi}m+\eta^1,\dots,\frac{\xi}m+\eta^{m-1},\frac{\xi}m-\sum_{j=1}^{m-1}\eta^j\big).
 \end{equation}
(In  the bilinear case ($m=2$), it is convenient to use the more symmetric parametrization of $ \rho_{H_0}(\xi,\eta)=(\frac{\xi-\eta}2,\frac{\xi+\eta}2)$.) One can then express the norm of $K$ in (\ref{L2new}) as

\begin{equation}\label{KrhoH}
||K||_{\Lambda^{H_0}_{2,\frac{2r}{r-2}}}=
c{\left[ \int_{\R^d} {\left( \int_{\R^{(m-1)d}} {\left|\widehat{K}\left(\rho_{H_0}\left(\xi,\veta\right)\right)\right|}^2 d\veta \right)}^{\frac{r}{r-2}} d\xi \right]}^{\frac{r-2}{2r}}. 
\end{equation}

Under the Fourier transform, Theorem \ref{bilinearnewwithattitude} becomes a bilinear Fourier restriction estimate. 
Interchanging the roles of $F,G$ and $\widehat{F},\widehat{G}$, consider $H$ as a subspace of $(\R^n)^*$ and, as with the diagonal $H_0\subset \R^{md}$ above,   introduce linear coordinates $u\in\R^k,\, v\in\R^{n-k}$ on $H^*,\, H^\perp$.
This gives rise to a linear isomorphism $\rho_H:\R^k_u\times \R^{n-k}_v\to\R^n$.  Thm. \ref{bilinearnewwithattitude}  then yields:

\vskip.125in 

\begin{corollary} \label{restriction} Let $H\subset (\R^n)^*$ be a $k$-dimensional subspace, $0\le k\le n$, and $F,G \in {\mathcal S}({\Bbb R}^{n})$. 
Then, for $p,q,r \ge 2$, $\frac{1}{p}+\frac{1}{q}+\frac{1}{r}=1$, we have 

$$ {|| \widehat{FG}|_H ||}_{L^r(H)} \lesssim {||F \circ \rho_H||}_{L^p_uL_v^2} \cdot {||G \circ \rho_H||}_{L_u^qL_v^2},$$
where $L_u^pL_v^2$ denotes the mixed-norm space with norm $\left(\int_{\R^k} \left(\int_{\R^{n-k}} |F(u,v)|^2 dv\right)^\frac{p}2 du\right)^\frac1{p}$.

\end{corollary}

\vskip.125in 

\begin{remark} Note that the extreme cases in Cor. \ref{restriction} of $dim(H)=n$ or $0$ correspond to standard facts: For  $H=\R^n$, one has $ H^\perp=(0)$ and $ L^p_uL^2_v=L^p_u$, so that  $FG\in L^{r'}$ by H\"older and then $\widehat{FG}\in L^r$ by Hausdorff-Young, while for  $H=(0),\, H^\perp=\R^n,\, L^p_uL^2_v=L^2_v$ and $FG\in L^1\implies \widehat{FG}(0)$ is well-defined.
\end{remark}

Starting with a single function $F$, replacing  $F$ and $G$ in Cor. \ref{restriction}  by $|F|^\frac12$ and $ |F|^\frac12\cdot\frac{F}{|F|}$, resp.,  convolving each with an approximate identity to restore membership in $\mathcal S$, and passing to the limit, one obtains a linear restriction theorem for linear subspaces:

\begin{corollary} \label{restrictiontohyperplanes} If  $H\subset (\R^n)^*$ is a $k$-dimensional subspace and $F \in {\mathcal S}({\Bbb R}^{n})$, 
then, for $p,q,r\ge 2,\, \frac1p+\frac1q+\frac1r=1$,

\begin{equation} \label{scramble} 
\big|\big|\widehat{F}|_{H}\big|\big|_{L^r(H)} \lesssim \big|\big|F\circ\rho_H\big|\big|^{\frac12}_{L^p_uL^1_{v}}\cdot \big|\big|F\circ\rho_H\big|\big|^{\frac12}_{L^q_uL^1_{v}},\quad p,q,r\ge 2,\, \frac1p+\frac1q+\frac1r=1,
\end{equation} 
and

\begin{equation} \label{linearrestriction} 
\big|\big|\widehat{F}|_{H}\big|\big|_{L^r(H)} \lesssim \big|\big|F\circ\rho_H\big|\big|_{L^p_uL^1_{v}},\quad r=\frac{p}{p-2},\, p\ge 2.
\end{equation}

\end{corollary} 

\vskip.125in 

Corollary \ref{restrictiontohyperplanes} shows that the Fourier transform of a  function can be restricted to a \emph{linear} subspace, provided that the mixed norm  is finite; this is in contrast with the well-known fact that the Fourier transforms of general  functions in  a standard $L^p,\, p>1,$ cannot be restricted to subspaces. 

\vskip.125in

%\vskip.125in 

Thm. \ref{bilinearnewwithattitude} also has a maximal operator variant: 

\begin{theorem} \label{maximalnew} Let $F,G \in {\mathcal S}({\Bbb R}^{n})$ and, for $t>0$,  define $G_t(\vx)=t^{-n} G\left( \frac{\vx}{t} \right)$. Let $H\subset {\Bbb R}^{n}$ be a $k$-dimensional subspace. Suppose that, for some $\gamma>\frac{1}{2}$  and any $t>0$, 

\begin{equation}\label{lpdecay}
||G^j_t||_{\Lambda^H_{2,\infty}}:= \sup_{\xi \in H^*} {\left( \int_{H_{\xi}^{\perp}} {\left|\widehat{G}^j(t\vzeta)\right|}^2  dH_{\xi}^{\perp}(\vzeta) \right)}^{\frac{1}{2}} \lesssim 2^{-j \gamma},
 \end{equation}
where, for $j\ge 0$,  $G^j$ is the Littlewood-Paley component of $G$ at frequency scale $2^j$ relative to a dyadic decomposition of $(\Bbb R^{n})^*$, and the same estimate holds for $\nabla\widehat{G}^j$. Then, 

\begin{equation} \label{maximalnewest} {\left|\left| \sup_{t>0} \left|F*G_t\right| \right|\right|}_{L^2(H)} \lesssim{||F||}_{L^2({\Bbb R}^{md})}.
\end{equation} 
\end{theorem} 

\vskip.125in 

\begin{remark} \label{garbageingarbageout} Approximating  $F$ by a continuous function plus a function with a small $L^2$ norm,
a  standard argument then shows that, for $F\in L^2(\R^{n})$, 

$$ \lim_{t \to 0} \int_{\R^n} F(x-u) G_t(u)\, du=c_G\cdot F(x)\, a.e.,\hbox{ with } c_G:= \int_{\R^{n}} G(u)\, du. $$  \end{remark} 

\vskip.125in 

Returning to the multilinear setting and 
taking $F(\vx)=\prod_{j=1}^m f_j(x^{j})$, $G(\vx)=K(\vx)$, $H=H_0$, and $p=2$ in Thm. \ref{maximalnew}, we obtain the following bound for maximal $m$-linear  convolution operators. 

\begin{corollary} \label{bilinearmaximalnew} If 
$K$ is a finite Borel measure on ${\Bbb R}^{md}$, 
$f_j \in {\mathcal S}({\Bbb R}^{d}),\, 1\le j\le m$, and $t>0$, define

$$B_t(f_1,\dots,f_m)(x)=\int_{\R^{md}} f_1(x-tu^{1})\cdots f_m(x-tu^{m})\, dK(\vu),$$ 
and the associated maximal operator,  
$$ {\mathcal B}(f_1,\dots,f_m)(x)=\sup_{t>0} |B_t(f_1,\dots,f_m)(x)|.$$ 
Then  

\begin{equation} \label{L2maximalnew} {||{\mathcal B}(f_1,\dots,f_m)||}_{L^2(\R^d)} \lesssim\prod_{j=1}^m {||f_j||}_{L^2(\R^d)}, \end{equation} 
provided that $K$ and $\nabla K$ satisfy the Littlewood-Paley   decay condition (\ref{lpdecay}) w.r.t. $\Lambda^{H_0}_{2,\infty}$.

\end{corollary} 

In Sec. \ref{maximalexamples} below, we use Cor. \ref{L2maximalnew} and some additional machinery to establish $m$-linear analogues of Stein's spherical maximal theorem \cite{St76}. 

\vskip.125in 

We point out  that there have been other works on restrictions of convolutions (e.g., \cite{BCW,BHT}) and estimates for multilinear convolution operators \cite{gressman,stovall}, but the results there seem to be of a different nature than those we obtain.

The paper is organized as follows. In Sec. \ref{bilinearapplications} we use Thm. \ref{bilinearnewwithattitude} to prove estimates for multilinear analogues of both Stein's spherical maximal theorem and  classical $L^p$-improving estimates for convolution.
The proofs    of Theorems \ref{bilinearnewwithattitude} and \ref{maximalnew} are then given in Sec. \ref{proofs sec}.
We  show in Sec. \ref{applications} how these results  imply certain Sobolev trace  inequalities and estimates for solutions of the heat and wave equations. 
 In Sec. \ref{oscillatorysection} we use our method to derive bounds for multilinear oscillatory integral operators, 
 and  Sec.   \ref{thoughts} contains some further discussion and open problems.

\vskip.225in 

\section{$L^p$ improving and maximal estimates for multilinear operators: }
\label{bilinearapplications}

\vskip.125in 

Multilinear operators play a key role in many aspects of harmonic analysis and partial differential equations. See, for example, \cite{St93} and the references there for a comprehensive description and history of this subject. The key starting point for us is the work of Coifman and Meyer \cite{CM90}, who initiated the comprehensive study of Calder\'on-Zygmund theory in the multilinear setting. They showed that an $m$-linear multiplier operator 
$$M: L^{p_1}({\Bbb R}^d) \times \dots \times L^{p_m}({\Bbb R}^d) \to L^r_s({\Bbb R}^d)$$ if 
$$ \frac{1}{p_1}+\dots+\frac{1}{p_m}=\frac{1}{r}$$ and  the multiplier $m(\vxi)$ of $M$ satisfies 
$$ |D^{\alpha}m(\xi^1, \dots, \xi^m)| \leq C_{\alpha} {(1+|\xi^1|+\dots+|\xi^m|)}^{-s-|\alpha|},\, \forall \alpha\in \Bbb Z_+^{md}.$$

This estimate has many important and interesting applications. However, in many multilinear situations  the symbolic decay conditions on the multiplier and its derivatives  are violated. 
For example, if the multiplier is the Fourier transform of the Lebesgue measure on a smooth hypersurface, the derivatives satisfy the same decay estimates as the multiplier, and no better. 
Another example which  is not covered by this result is the case of multilinear fractional integration, handled in \cite{KSt99} using different methods; see also \cite{GK01}. 

A systematic study of multilinear analogues of linear generalized Radon transforms was begun in \cite{GGIP12}, focussing on applications to the study of  finite point configurations problems in geometric measure theory. These problems have their roots in the work of Falconer on the distance problem in the continuous setting and of the Erd\H os school on finite point configurations in the discrete setting. See, e.g., \cite{BMP05} and the references there. The main bilinear estimate in \cite{GGIP12} can be stated as follows. 

\begin{theorem} \label{bilineargrandma} Let $K$ be a {\it non-negative} integrable function on ${\Bbb R}^{2d} $, and define

$$ B(f,g)(x)=\int f(x-u) g(x-v) K(u,v) du dv. $$ 
Then for $1 \leq p \leq 2$, 
$$ {||B(f,g)||}_{L^p({\Bbb R}^d)} \leq {||f||}_{L^2({\Bbb R}^d)} \cdot {||g||}_{L^2({\Bbb R}^d)} \cdot {\left( \int {\left| \widehat{K}(\xi, -\xi) \right|}^{p'} d\xi \right)}^{\frac{1}{p'}}. $$ 
\end{theorem} 

While this result was useful in establishing geometric and combinatorial estimates, a weakness that limits wider applications, is the positivity assumption on the kernel. This prevents one  from exploiting, e.g., Fourier decay information by using Littlewood-Paley type decompositions. We resolve this problem to a significant extent in Thm. \ref{bilinearnewwithattitude} and Cor. \ref{bilinearnew},   opening the door to a variety of applications illustrated in this section. 

\subsection{$L^p$-improving measures} 
\label{lpimproving}

A finite measure $\mu$ on ${\Bbb R}^d$ is said to be \emph{$L^p$-improving} if 
\begin{equation} \label{Lpimprovingexponent} {||f*\mu||}_{L^q({\Bbb R}^d)} \leq C_{p,q} {||f||}_{L^p({\Bbb R}^d)}\hbox{ for some } q>p. \end{equation} 
  
There has been much interest in determining the pairs $(p,q)$ for which (\ref{Lpimprovingexponent}) holds.
Two representative results in this area are due to Strichartz \cite{Str70} and Littman \cite{L71}; see Christ \cite{C98} for more recent work. The result of Littman and Strichartz says that if $\mu$ is the Lebesgue measure on the unit sphere in ${\Bbb R}^d$, $d \ge 2$, then (\ref{Lpimprovingexponent}) holds if and only if $(1/p, 1/q)$ is contained in the closed triangle with endpoints $(0,0)$, $(1,1)$ and $\left(\frac{d}{d+1}, \frac{1}{d+1} \right)$. 

\vskip.125in 

In the multilinear setting, the notion of $L^p$-improving may be expressed in the following way, indicating an improvement over the H\"older inequality exponents governing pointwise multiplication.

\begin{definition} Let $B(f_1,\dots,f_m)$ be an $m$-linear operator. We say that $B$ is {\it Lebesgue space improving} if there exist $p_1,\dots,p_m,r>0$ with 

$$ \frac{1}{p_1}+\cdots +\frac{1}{p_m}>\frac{1}{r}$$ such that 

$$ B: L^{p_1}({\Bbb R}^d) \times\cdots \times L^{p_m}({\Bbb R}^d) \to L^r({\Bbb R}^d).$$ 
\end{definition} 

\vskip.125in 

The following result is a simple adaptation of Cor. \ref{bilinearnew}. 

\begin{theorem} \label{Lpimproving} Suppose that $\nu(\vu)$ is a smooth multiple of surface measure on a compact, codimension $l$ submanifold of ${\Bbb R}^{md} $. For $j\ge 0$, let $\nu_j$ be the $j$th Littlewood-Paley component of $\nu$, and suppose that, for some $\gamma>0$, 

\begin{equation} \label{kaif} 
||\nu_j||_{\Lambda^{H_0}_{2,\infty}}\lesssim 2^{-j\gamma}. \end{equation} 
Define

$$B_{\nu}(f_1,\dots,f_m)(x) := \int_{\R^{md}} f_1(x-u^1)\cdots f_m(x-u^m)\, d\nu(\vu).$$ 
Then 

\begin{equation} \label{sobolevgeneralestimate} 
{||B_{\nu}(f_1,\dots,f_m)||}_{L^2_{\gamma}({\Bbb R}^d)} \lesssim \prod_{j=1}^m {||f_j||}_{L^2({\Bbb R}^d)}  \end{equation} 
and

\begin{equation} \label{Lpimprovinggeneralestimate}
 {||B_{\nu}(f_1,\dots,f_m)||}_{L^{p'}({\Bbb R}^d)} \lesssim \prod_{j=1}^m {||f_j||}_{L^p({\Bbb R}^d)} 
\ \text{for} \ p>\frac{2(l+\gamma)}{l+2 \gamma}. \end{equation} 
\end{theorem} 

\vskip.125in 

\begin{remark} \label{pointwisetoaverage} It is not difficult to check that if $|\widehat{\nu}(\vxi)| \leq C{(1+|\vxi|)}^{-\alpha}$, then (\ref{kaif}) holds with the exponent $\gamma=\alpha-\frac{(m-1)d}{2}$. \end{remark} 

\begin{corollary} \label{fullsphere} Let $\nu$ be  surface measure on the sphere $\mathbb S^{md-1}=\big\{\vu \in {\Bbb R}^{md} : {|u^1|}^2+\cdots+{|u^m|}^2=\nolinebreak1\big\}$. Then 

$$ B_{\nu}: L^p({\Bbb R}^d) \times\cdots\times L^p({\Bbb R}^d) \to L^{p'}({\Bbb R}^d) \ \text{for} \ \frac{d+1}{d}<p \leq 2.$$ 
\end{corollary} 

\vskip.125in 

Cor. \ref{fullsphere} follows from Thm. \ref{Lpimproving}, since the codimension $l=1$ and the estimate (\ref{kaif}) holds with $\gamma=\frac{md-1}{2}-\frac{d}{2}$, using the Remark \ref{pointwisetoaverage} and the 
well-known decay (see, e.g., \cite{St93}) 
that, if  $\sigma$ denotes surface measure on the sphere in ${\Bbb R}^n$, $n \ge 2$, then 
\begin{equation}\label{spheredecay}
 |\widehat{\sigma}(\xi)| \leq C{(1+|\xi|)}^{-\frac{n-1}{2}}.
 \end{equation} 

One can in fact improve Cor. \ref{fullsphere} to obtain the endpoint  $p=\frac{d+1}d$ as follows. Define 

$$ B_{\nu}^z(f_1,\dots,f_m)(x):=\int  f_1(x-u^1) \cdots f_m(x-u^m)\,  d\nu^z(\vu),$$ where 
$$ \nu^z(\vu)=\frac{1}{\Gamma(z)} {(1-{|\vu|}^2)}_{+}^{z-1}.$$ 

\vskip.125in

When $Re(z)=-\frac{d-1}{2}$, $B^z_{\nu}: L^2({\Bbb R}^d) \times\cdots\times L^2({\Bbb R}^d) \to L^2({\Bbb R}^d)$ via a direct appeal to Cor. \ref{bilinearnew} and the observation above that we may take all the way up to $\gamma=\frac{(m-1)d-1}{2}$ in (\ref{kaif}). When $Re(z)=1$, $B_{\nu}: L^1({\Bbb R}^d) \times\cdots\times L^1({\Bbb R}^d) \to L^{\infty}({\Bbb R}^d)$ trivially since the kernel becomes bounded. Stein's analytic interpolation theorem then yields the conclusion of Cor. \ref{fullsphere} without the loss of the endpoint.

This estimate is basically sharp, as can be seen by the following adaptation of the usual linear sharpness example which we illustrate in the bilinear setting. Let $f$ be the characteristic function of the unit ball and let $g$ be the characteristic function of the ball of radius $\delta$. Then ${||g||}_p \approx \delta^{\frac{d}{p}}$ and $B_{\nu}(f,g) \approx \delta^{d-1}$ on the annulus of thickness $ \approx \delta$ and radius $\approx 1$. Sharpness of the estimate above, up to the endpoint, follows.

\vskip.125in 

In the bilinear case, if $\{(u^1,u^2): {|u^1|}^2+{|u^2|}^2=1 \}$ is replaced by $\{(u^1,u^2): u^1 \cdot u^2=t \}$ with $t \not=0$ fixed, one obtains
the same exponents as in Cor. \ref{fullsphere}. This surface arose in \cite{GGIP12} in connection with the study of areas of triangles determined by fractal subsets of  Euclidean space. 

\vskip.125in 

\subsection{Proof of Theorem \ref{Lpimproving}} 

\vskip.125in

The estimate (\ref{sobolevgeneralestimate})  follows  from Thm. \ref{bilinearnewwithattitude} and (\ref{kaif}). 
In fact, since $B_\nu(f_1,\dots,f_m)=(f_1\otimes\cdots\otimes f_m)*\nu|_{H_0}$,  we have $(I-\Delta_x)^{\gamma/2} B_\nu(f_1,\dots,f_m)=
(f_1\otimes\cdots\otimes f_m)*V|_{H_0}$, where $\widehat{V}(\xi^1,\dots,\xi^m)=(1+|\xi^1+\cdots+\xi^m|^2)^{\gamma/2}\widehat{\nu}$.
From (\ref{kaif}), we see that $||V||_{\Lambda^{H_0}_{2,\infty}}\lesssim 1$, and hence applying Thm. \ref{bilinearnewwithattitude} with $r=p=2,\, q=\infty$ yields (\ref{sobolevgeneralestimate}).

\vskip.125in

To prove (\ref{Lpimprovinggeneralestimate}), define as above the $j$th Littlewood-Paley piece $\nu_j$ of $\nu$,  $0\le j<\infty$, and set 
$$B^j(f_1,\dots,f_m)(x)=\int\cdots  \int f_1(x-u^1)\cdots  f_m(x-u^m)\, \nu_j(\vu)d\vu.$$
It follows from (\ref{sobolevgeneralestimate}) that 

$$ B^j: L^2({\Bbb R}^d) \times\cdots\times L^2({\Bbb R}^d) \to L^2({\Bbb R}^d) \ \text{with norm} \ 2^{-j \gamma}.$$ 
On the other hand, 

$$ {||B^j(f_1,\dots,f_m)||}_{L^{\infty}({\Bbb R}^d)} \leq {||f_1||}_{L^1({\Bbb R}^d)} \cdots {||f_m||}_{L^1({\Bbb R}^d)} \cdot {||\nu_j||}_{L^{\infty}({\Bbb R}^{md}) }.$$ 
By a direct calculation (see, e.g., \cite{W03}), using the codimension $l$ of the support of $\nu$,

$$ {||\nu_j||}_{L^{\infty}({\Bbb R}^{md})} \leq C2^{lj}.$$
The Riesz-Thorin interpolation theorem then yields (\ref{Lpimprovinggeneralestimate}).

\begin{corollary} \label{productsphere} Let $\nu$ be surface measure on $S^{d-1} \times\cdots \times S^{d-1} \subset {\Bbb R}^{md} $. Then 
$$ B_{\nu}: L^p({\Bbb R}^d) \times\cdots\times  L^p({\Bbb R}^d) \to L^{p'}({\Bbb R}^d) \ \text{for} \ \frac{(m^2-2m+2)d+m^2}{(m^2-2m+2)d} \leq p \leq 2.
$$ \end{corollary} 

To prove this, define $\nu^{z}(\vu)=\prod_{j=1}^m\frac{1}{\Gamma(z)} {(1-{|u^j|}^2)}_{+}^{z-1}$ and proceed as in the proof of Cor. \ref{fullsphere}.

\vskip.25in 

\subsection{Bilinear analogues of Stein's spherical maximal theorem} 
\label{maximalexamples}

The classical (linear)  spherical maximal operator is given by 

$$ {\mathcal M}f(x)=\sup_{t>0} \left| \int_{S^{d-1}} f(x-ty) d\sigma(y) \right|,$$ where $\sigma$ is th Lebesgue measure on the unit sphere in ${\Bbb R}^d$, $d \ge 2$. It was established by Stein \cite{St76} in three dimensions and higher, and by Bourgain \cite{B86} in two dimensions that  

$$ {\mathcal M}: L^p({\Bbb R}^d) \to L^p({\Bbb R}^d) \ \text{for} \ p>\frac{d}{d-1}.$$ 
The result is best possible, as can be shown by taking 

$$f(x)={|x|}^{-d+1} \log \left(\frac{1}{|x|} \right) \chi_{B(\vec{0}, \frac{1}{2})}(x) \ \text{with} \ t=|x|.$$ 

\vskip.125in

The cornerstone of the general theory of maximal averages in dimension $d\ge 3$ is the following result; see 
Bourgain \cite{B86II}, Carbery \cite{C86}, Cowling-Mauceri \cite{CM86} and Sogge-Stein \cite{SoSt85}. 
\begin{theorem} \label{maximallineargeneraltheorem} Let  $\mu$ be a finite, compactly supported measure on ${\Bbb R}^d$, such that, for some 
$\epsilon>0$, $|\widehat{\mu}(\xi)| \leq C{(1+|\xi|)}^{-\frac{1}{2}-\epsilon}$. Define

\begin{equation} \label{maximallineargeneral} {\mathcal M}_{\mu}f(x)=\sup_{t>0} \left| \int f(x-ty) d\mu(y) \right|. \end{equation} 
Then,  $ {\mathcal M}_{\mu}: L^2({\Bbb R}^d) \to L^2({\Bbb R}^d).$ 
\end{theorem} 

\vskip.125in

The techniques  introduced above  yield multilinear analogues of these maximal theorems. For simplicity, we restrict ourselves to the bilinear case.
Given a finite Borel measure $\nu$ on $\Bbb R^{2d}$, define the bilinear averaging operators,   

$$ B_t(f,g)(x)=\int \int f(x-tu) g(x-tv) d\nu(u,v),\, 0<t<\infty,$$ 
and the associated maximal operator, $ {\mathcal B}_{\nu}(f,g)(x)=\sup_{t>0} |B_t(f,g)(x)|$. 
The following bilinear analogue of Thm. \ref{maximallineargeneraltheorem} is an immediate restatement of Cor. \ref{bilinearmaximalnew}.

\begin{theorem} \label{bilinearftc} Suppose that, for all $t>0$ and some $\epsilon>0$,

\begin{equation} \label{bilinearaveragedecay} 
\sup_{\xi\in\Bbb R^d} {\left( \int_{2^j \leq \sqrt{{|t\xi|}^2+{|t\eta|}^2} \leq 2^{j+1}} 
{\left|\, \widehat{\nu}\left(\frac{t\xi+t\eta}{2}, \frac{t\xi-t\eta}{2}\right)\right|}^2 d\eta \right)}^{\frac{1}{2}} \leq C 2^{-j\left( \frac{1}{2}+\epsilon \right)}, \end{equation}
and the same estimate holds for $|\nabla\hat\nu|$. Then $ {\mathcal B}_{\nu}: L^2({\Bbb R}^d) \times L^2({\Bbb R}^d) \to L^2({\Bbb R}^d)$. 
\end{theorem} 

Using  additional geometric information on the support of $\nu$, one can improve the exponents in Thm. \ref{bilinearftc} as follows. 

\begin{corollary} \label{interpolation} Suppose that $\nu$ is supported on a smooth compact submanifold of ${\Bbb R}^d \times {\Bbb R}^d$ of codimension $l$. Suppose that, for all $t>0$ and some $\gamma>\frac{1}{2}$, 

\begin{equation} \label{bilinearaveragedecayconcrete} \sup_{\xi\in\Bbb R^d} {\left( \int_{2^j \leq \sqrt{{|t\xi|}^2+{|t\eta|}^2} \leq 2^{j+1}} 
{\left|\,\widehat{\nu}\left(\frac{t\xi+t\eta}{2}, \frac{t\xi-t\eta}{2}\right)\right|}^2 d\eta \right)}^{\frac{1}{2}} \leq C 2^{-j \gamma},\, j \ge 0 \end{equation}  
and the same holds for the gradient. Then 

\begin{equation} \label{maximalexponent} {\mathcal B}_{\nu}: L^p({\Bbb R}^d) \times L^p({\Bbb R}^d) \to L^{p'}({\Bbb R}^d) \ \text{with} \ p>\frac{2l+2\gamma-1}{l+2 \gamma-1}. \end{equation}  

\end{corollary} 

\vskip.125in 

\begin{remark} By a standard argument involving the approximation of an $L^p$ function by a continuous function plus a function with a small $L^p$ norm, Thm. \ref{bilinearftc} and Cor. \ref{interpolation} yield differentiation theorems. For example, Cor. \ref{interpolation} implies that if $f,g \in L^p({\Bbb R}^d)$ with $p$ as in  (\ref{maximalexponent}), then 
$$ \lim_{t \to 0} \int \int f(x-tu) g(x-tv) d\nu(u,v)=f(x)g(x) \ a.e. $$ 
\end{remark} 

\vskip.125in 

\begin{example} \label{fullspheremaximal} Let $\nu$ be the Lebesgue measure on the sphere $\{(u,v) \in {\Bbb R}^d \times {\Bbb R}^d: {|u|}^2+{|v|}^2=1\}$, $d \ge 2$. Then 

$$ {\mathcal B}_{\nu}: L^p({\Bbb R}^d) \times L^p({\Bbb R}^d) \to L^{p'}({\Bbb R}^d) \ \text{for} \ p>\frac{d}{d-1}.$$ 
\end{example} 

To see that the exponent $\frac{d}{d-1}$ is optimal, take $g \equiv 1$. What we get is essentially the maximal averaging operator of $f$ over the sphere in ${\Bbb R}^d$ of radius $1-{|v|}^2$. The needed restriction follows from Stein's sharpness example for the linear spherical maximal operator, namely, $f(x)={|x|}^{-d+1} \log \left(\frac{1}{|x|} \right) \chi_{B(\vec{0}, \frac{1}{2})}(x)$. 

\subsection{Proof of Corollary \ref{interpolation}} 

\vskip.125in 

By a direct calculation, 

$$ |B_t^j(f,g)(x)| \leq {||f||}_1 \cdot {||g||}_1 \cdot {||K^j||}_{\infty} \leq C {||f||}_1 \cdot {||g||}_1 \cdot 2^{jl}.$$

\vskip.125in 

The proof of Thm. \ref{bilinearftc} above gives 
$$ {\mathcal B}^j: L^2({\Bbb R}^d) \times L^2({\Bbb R}^d) \to L^2({\Bbb R}^d) \ 
\text{with norm} \ \approx 2^{-j \left(\gamma-\frac{1}{2} \right)},$$
and using Riesz-Thorin again yields Cor. \ref{interpolation}.

\vskip.25in

\section{Proofs    of Theorems \ref{bilinearnewwithattitude} and \ref{maximalnew}}\label{proofs sec}

\vskip.125in 

\subsection{Proof    of Thm. \ref{bilinearnewwithattitude}}

By the translation-invariance of the norms in (\ref{newyoung}), one may assume that $H\subset\R^{n}$ is a $k$-dimensional linear subspace, which,  by rotation covariance, can be put in the form $H=\{(x',x'')\in\R^k\times \R^{n-k} \,|\, x''=0\}\sim\R^k$, so that, for $\xi'\in H^*$, $H_\xi^\perp=\{(\xi',\xi'')\, |\, \xi''\in \R^{n-k}\}$. Fixing $G$, define the restricted convolution operator $TF=(F*G)|_H$, i.e., $TF(x')=(F*G)(x',0)$. Then the formal adjoint is given by

$$T^*h(u',u'')=\int_{\R^k} \overline{G}(x'-u',-u'')\, h(x')\, dx',$$
from which one sees that $TT^*$ is a translation-invariant operator on $\R^k$, given in Fourier multiplier form by

\begin{eqnarray*}
\widehat{TT^{*}h}(\xi') =  \left( \int_{\R^{n-k}}|\widehat{G}(\xi',\xi'')|^2 d\xi'' \right)\widehat{h}(\xi')= \left( \int_{H_{\xi}^{\perp}}|\widehat{G}(\vzeta)|^2 dH^\perp_{\xi'}(\vzeta) \right)\widehat{h}(\xi).
\end{eqnarray*}
By Parseval,

$$||TT^*||_{L^2(H)\to L^2(H)}=\sup_{\xi'\in H^*} \left( \int_{H_{\xi}^{\perp}}|\widehat{G}(\vzeta)|^2 dH^\perp_{\xi'}(\vzeta) \right) = ||G||^2_{\Lambda^H_{2,\infty}},$$
and thus $||T||_{L^2(H)\to L^2(\R^n)} = ||G||_{\Lambda^H_{2,\infty}}$, yielding (\ref{bestconstant}), and therefore  (\ref{bilinearnewwithattitude})  for $p=2,\, q=\infty,\, r=2$. On the other hand, 

$$ {||F*G||}_{L^{\infty}(H)} \le ||F*G||_{C_0(\R^{n})} \le ||\hat{F}\hat{G}||_{L^1(\R^{n})}\le ||\hat{F}||_{L^2}\cdot ||\hat{G}||_{L^2}
\leq {||F||}_{\Lambda^H_{2,2}} \cdot {||G||}_{\Lambda^H_{2,2}},$$ 
which is (\ref{bilinearnewwithattitude}) for $p=q=2,\, r=\infty$.
Interpolation then gives (\ref{bilinearnewwithattitude}) in for general $p,q,r$.

\vskip.125in 

\subsection{Proof of Thm. \ref{maximalnew}}

It is sufficient to establish
\begin{equation}\label{supsup}
\| \sup\limits_R \ \sup\limits_{R\leq t \leq 2R} |F*G_t| \|_{L^2(H)}^{2} \leq \| F \|^2_{2}
\end{equation}
where $R$ ranges over all dyadic numbers. For fixed $R$, make a Littlewood-Paley decomposition with respect to that scale. First define 
\begin{equation} \label{littlewoodpaley} 
\widehat{F^j}(\vxi)=\psi\left(2^{-j} |\vxi|\right)\widehat{F}(\vxi), \end{equation} 
where $\psi$ is a smooth cut-off function supported in the annulus in ${\Bbb R}^{md} $ of inner radius $\frac{1}{2}$ and outer radius $4$, identically equal to $1$ in the annulus of inner radius $1$ and outer radius $2$, and satisfying $\sum_j \psi(2^{-j} \cdot) \equiv 1$. Letting $[\cdot]$ denote the greatest integer function, we now aim to estimate
$$  \sup\limits_{R\leq t \leq 2R} |F^{j+[\log_2(R)]}*G_t| .$$

We shall need the following elementary observation. 
\begin{lemma} \label{ftc} Let $F$ be a differentiable function on $[R,2R]$. Then 
$$ \sup_{t \in [R,2R]} {|F(t)|}^2 \leq {|F(R)|}^2+2 {\left( \int_{R}^{2R} {|F(t)|}^2 dt \right)}^{\frac{1}{2}} \cdot {\left( \int_R^{2R} {|F'(t)|}^2 dt \right)}^{\frac{1}{2}}.$$ 
\end{lemma} 

To prove the lemma, just observe that by the fundamental theorem of calculus, 
$$ {F(t)}^2={F(R)}^2+2 \int_R^t 2F(s)F'(s)ds$$ and apply the Cauchy-Schwarz inequality. 
In order to use Lemma \ref{ftc}, observe that 
$$ \frac{d}{dt} (F * G_t)(x)=\int e^{2 \pi i x \cdot (\xi^1+\cdots+\xi^m)} \widehat{F}(\vxi) \frac{d}{dt} 
\widehat{G}(t\vxi) d\vxi $$
$$=\int e^{2 \pi i x \cdot (\xi^1+\cdots+\xi^m)} \widehat{F}(\vxi) \langle\nabla \widehat{G}(t\vxi), \vxi\rangle d\vxi,$$
$$= \int e^{2 \pi i x \cdot (\xi^1+\cdots+\xi^m)} \widehat{F}(\vxi)t^{-1}\widehat{G^{*}}(t\vxi) d\vxi ,$$ where 
$$ \left|\widehat{G^{*}}(\vxi)\right| \lesssim |\vxi| \cdot \left|\nabla\widehat{G}(\vxi)\right|.$$ 

Note that
\begin{align*}
\left(F^{j+[\log_2(R)]}*G_t \right)(x) &=  \int e^{2 \pi i x \cdot (\xi^1+\cdots+\xi^m)} \widehat{F}(\vxi)\psi\left(2^{-j} R |\vxi|\right) \widehat{G}(t\vxi) d\vxi  \\
&= \int e^{2 \pi i x \cdot (\xi^1+\cdots+\xi^m)} \widehat{F}(\vxi)\psi\left(2^{-j} t |\vxi|\right) \widehat{G}(t\vxi) d\vxi  \\
&= \left(F*G_t^j \right)(x)
\end{align*}
when $R\leq t \leq 2R$, where  $G_t^j$ is given by
$ \widehat{G_t^j}(\vxi) = \psi\left(2^{-j} t |\xi|\right) \widehat{G}(t\vxi)$.
Similar calculations hold for \nolinebreak$G^*$.

\vskip.125in 

Applying Lemma \ref{ftc}, we see that 

\begin{eqnarray*}
 {||\sup\limits_{t \in [R,2R]}|F^{j+[\log_2(R)]}*G_t|_H| ||}^2_{L^2(H)} &\leq&  {||F*G^j_R|_H||}_{L^2(H)}^2\\
& &+2 {\left(\frac{1}{R} \int_R^{2R} \int {|(F*G_t^j)(x)|}^2 dxdt \right)}^{\frac{1}{2}}\\
& \quad  & \times
{\left(R \int_R^{2R}t^{-2} \int {|(F*(G^{*})^j_t)(x)|}^2 dxdt \right)}^{\frac{1}{2}} .
\end{eqnarray*}

\vskip.125in 

Now appealing to  Thm. \ref{bilinearnewwithattitude} and using estimate (\ref{lpdecay}), we see that this expression is 

$$ \lesssim {||F||}_2^2 \cdot 2^{-j(1+2 \epsilon)} +
{||F||}_2^2 \cdot 2^{-j \left( \frac{1}{2}+\epsilon \right)} \cdot 2^{-j \left( \frac{1}{2}+\epsilon \right)} 2^j.$$ 

\vskip.125in 

We could have made the argument above stronger by replacing the $\psi$ factor with
$ \psi\left(2^{-j} R |\xi|\right)\tilde{\psi}\left(2^{-j} R |\xi|\right)$,
where $\tilde{\psi}$ is a slight widening of $\psi$. The  argument as above would then yield
$$ {||\sup\limits_{t \in [R,2R]}|F^{j+[\log_2(R)]}*G_t|_H\, \|}^2_{L^2(H)} \lesssim 2^{-j\epsilon}{||F^{j+[\log_2(R)]}||}_2^2, $$
where the Littlewood-Paley decomposition on the right hand side is made with respect to the dilates of $\tilde{\psi}$.
Summing the corresponding geometric series and bounding a supremum by a square function, we obtain equation (\ref{supsup}) and thus the conclusion of Theorem \ref{maximalnew}. 

\vskip.25in

%%%

\section{Applications to PDE}\label{applications}

We now show that the results established so far can be used to obtain a variety of estimates related to Sobolev spaces and partial differential equations.

\subsection{Sobolev traces on subspaces}\label{Sobolev}

%\vskip.25in 

The classical Sobolev trace inequality from $\R^n$ to $\R^{n-1}$ (see, e.g.,  \cite[Prop. 1.6]{T96}) says that if the restriction (or trace) map $\tau$ is defined by $\tau(u)=f$, where $f(x')=u(0,x')$, $x=(x_1,x')$, $x'=(x_2, \dots, x_n)$, then 
\begin{equation} \label{taylor} {||\tau (u)||}_{L^2_{s-\frac{1}{2}}({\Bbb R}^{n-1})} \leq C_s{||u||}_{L^2_s({\Bbb R}^n)},\quad s>\frac{1}{2}. \end{equation}
Iterating (\ref{taylor}) $n-k$ times yields, for any $k$-dimensional subspace $H\subset \Bbb R^{n}$,
\begin{equation} \label{taylorsharp} {||\tau(u)||}_{L^2_{s-\frac{n-k}{2}}(H)} \leq C_s' {||u||}_{L^2_s({\Bbb R}^{2d})}, \quad s>\frac{n-k}2. \end{equation} 

We shall use Thm. \ref{bilinearnewwithattitude} to prove an endpoint version of (\ref{taylorsharp}) with the best constant, 
and then give some applications to
partial differential equations. 

\begin{theorem} \label{bestconstantresult} 
If $H\subset{\Bbb R}^{n} $  is an affine subspace with $dim(H)=k$, and $s>\frac{n-k}{2}$, then 

\begin{equation} \label{taylorsharpest} {||u|_{H}||}_{L^2(H)} \leq C_{s,H} {||u||}_{L^2_s({\Bbb R}^{n})},\, \forall u \in {\mathcal S}({\Bbb R}^{n}). \end{equation} 
Furthermore,  the optimal constant is

\begin{equation} \label{firstbestconstant} C_{s,n-k}=\big|\big|(1+|\xi|^2)^{-\frac{s}2}\big|\big|_{\Lambda^H_{2,\infty}}=\sqrt{\big|\mathbb S^{n-k-1}\big|\cdot  \int_0^{\infty} {(1+r^2)}^{-s} r^{n-k-1} dr}. \end{equation}  

\end{theorem} 

\vskip.125in

\begin{remark} Thm. \ref{bestconstantresult} is result is an endpoint version of (\ref{taylorsharp}) in the sense that we take $L^2(H)$ on the left hand side instead of $L^2_{s-\frac{n-k}{2}}(H)$ for $s>\frac{n-k}{2}$. One cannot pin the endpoint on both sides, as boundedness $\tau:L^2_{\frac{d}{2}}({\Bbb R}^{2d}) \to L^2(H)$  does not hold. Clearly, 
if  $C_{s,n-k}$ is as in (\ref{firstbestconstant}) then
$lim_{s \searrow \frac{n-k}{2}} \, C_{s,n-k}=\infty$. 
\end{remark} 

\vskip.125in

\noindent{\bf Proof of Theorem \ref{bestconstantresult}.} 
As in the proof of Thm. \ref{bilinearnewwithattitude}, by translation- and rotation-invariance of the function spaces involved, we can assume that $H=\{(x',x'')\in \R^k\times\R^{n-k}\, |\, x''=0\}$, so that $H^*\simeq \R^k_{\xi'}$ and $H^\perp_{\xi'}=\{(\xi',\xi'')\, |\, \xi''\in\R^{n-k}\}$.
One can rewrite (\ref{taylorsharpest}) in the form 
$$ {|| {(-\bigtriangleup)}^{-\frac{s}{2}}u ||}_{L^2(H)} \leq C_s {||u||}_{L^2({\Bbb R}^{2d})},\, s>\frac{d}{2}.$$
The left hand side can be rewritten as 
$$ {|| u*G_s ||}_{L^2(H)} \ \text{with} \ \widehat{G}_s(\xi',\xi'')={(1+{|\xi'|}^2+{|\xi''|}^2)}^{-\frac{s}{2}}.$$ 
By Thm. \ref{bilinearnewwithattitude} we conclude that  (\ref{taylorsharpest}) holds, 
with optimal constant 
\begin{eqnarray*}
C_{s,n-k}&=& \big|\big|(1+|\xi|^2+|\xi''|^2)^{-\frac{s}2}\big|\big|_{\Lambda^H_{2,\infty}}.  
\end{eqnarray*}
The $\sup_{\xi'}$ in the norm on the right is attained at $\xi'=0$, yielding

$$C_{s,n-k}^2= \big|\mathbb S^{n-k-1}\big|\cdot
\int_0^{\infty} (1+r^2)^{-s} r^{n-k-1}\, dr<\infty,\quad  s>\frac{n-k}2,$$
completing  the proof of (\ref{taylorsharpest}) with constant (\ref{firstbestconstant}). 

\vskip.125in 

We shall also give a simple proof of the following product Sobolev estimate \cite{T01}: 

\begin{theorem} \label{laceysobolevproduct} Let $u,v \in {\mathcal S}({\Bbb R}^d)$. Suppose that $\gamma=r+s-\frac{d}{2}$ and $r,s \ge 0$. Then 

\begin{equation}  \label{laceysobolevproductest}{||uv||}_{L^2_{\gamma}({\Bbb R}^d)} \leq C {||u||}_{L^2_r({\Bbb R}^d)} \cdot {||v||}_{L^2_s({\Bbb R}^d)}. \end{equation} 
\end{theorem} 

\begin{remark} Unfortunately, we are not able to obtain the best constant in (\ref{laceysobolevproductest}), since we cannot be certain that the optimizing function for the inequality (\ref{bestconstant}) is realizable as a product function. \end{remark} 

%%%%

\vskip.125in

\noindent{\bf Proof of Theorem \ref{laceysobolevproduct}.} 
We make use of the bilinear convolution estimates  from Cor. \ref{bilinearnew}. Let 
$$f={(-\bigtriangleup)}^{\frac{s}{2}}u \ \text{and} \ g={(-\bigtriangleup)}^{\frac{r}{2}}.$$ 
We can rewrite (\ref{laceysobolevproductest}) in the form 

$$ {||uv||}_{L^2_{\gamma}({\Bbb R}^d)} \leq C {||f||}_{L^2({\Bbb R}^d)} \cdot {||g||}_{L^2({\Bbb R}^d)}.$$ 
Let 
$$ B(f,g)(x)=c_{r,s} \int \int f(x-u) g(x-v) {|u|}^{-d+r} {|v|}^{-d+s} du dv$$ and observe that, with the appropriate choice of $c_{r,s}$, one has $B(f,g)=uv$. 
Applying Cor. \ref{bilinearnew}, we see that 
$$ {||uv||}_{L^2({\Bbb R}^d)} \leq {||f||}_{L^2({\Bbb R}^d)} \cdot {||g||}_{L^2({\Bbb R}^d)} \cdot 
\sup_{\xi \in {\Bbb R}^d} 
{\left( \int {\left|\widehat{K}\left(\frac{\xi+\eta}{2}, \frac{\xi-\eta}{2} \right)\right|}^2 d\eta \right)}^{\frac{1}{2}},$$ 
where  $ \widehat{K}(\xi, \eta)={|\xi|}^{-r} {|\eta|}^{-s}$.
Since 
$$  {\left( \int_{|\xi|+|\eta| \approx 2^j} {\left|\widehat{K}\left(\frac{\xi+\eta}{2}, \frac{\xi-\eta}{2} \right)\right|}^2 d\eta \right)}^{\frac{1}{2}} \approx 2^{-j(r+s-\frac{d}{2})},$$ the result follows.

\vskip.125in 

We now give some applications of Thms. \ref{bestconstantresult} and \ref{laceysobolevproduct} to  {\it a priori} estimates in partial differential equations.

\vskip.125in 

\subsection{The heat equation: restriction to subspaces} Consider the  heat equation on ${\Bbb R}^{n}\times\nolinebreak\mathbb R_+$,
\begin{equation} \label{heatequation} 
\frac{\partial u}{\partial t}=\bigtriangleup u, \, (x,t)\in\R^{n}\times\R_+,\, \ u(x,0)=F(x)\in L^2(\mathbb R^{n}),\end{equation} 
with  solution  given by
$ u(x,t)=\big(F*\Phi_t\big)(x)$, where 
$$ \Phi_t(x)={(4 \pi t)}^{-\frac{n}2} e^{-\frac{{|x|}^2}{4t}} .$$ 
We have the following result concerning the restriction of the heat semi-group to affine subspaces.  

\begin{theorem} \label{heatcontraction} Suppose that $u$ is the solution to  (\ref{heatequation}). Then, for any $k$-dimensional affine subspace $H\subset \mathbb R^{n}$,
$$ {\big|\big|u(\cdot,  t)|_H\big|\big|}_{L^2(H)} \leq {\left( {4 \pi t} \right)}^{-\frac{n-k}{4}} \cdot {\big|\big|F\big|\big|}_{L^2({\Bbb R}^{n})},$$ 
and the constant ${\left( {4 \pi t} \right)}^{-\frac{n-k}{4}}$ is optimal. 
\end{theorem} 

In other words,  the initial data $F$  might  blowup along  $H$ for $t$ small, resulting in the solution $u$ being large along $H$ for small values of $t$, but after the fixed time $T= \frac{1}{4\pi}$, independent of $F$, the $L^2(H)$ norm of the solution drops below the $L^2({\Bbb R}^{n})$ norm of $F$. 

\vskip.125in

\noindent{\bf Proof of Theorem \ref{heatcontraction}.} By invariance of (\ref{heatequation}) under the Euclidian motion group, as in the proof of Thm.  \ref{bestconstantresult}  $H$ may   be taken to be $H=\{(x',0)\in\R^n\, |\, x'\in\R^k\}$. Since 
$ \widehat{\Phi}_t(\xi)=e^{-4 \pi^2 t{|\xi|}^2}=e^{-4 \pi^2 t({|\xi'|}^2+{|\xi''|}^2)}$, one has

\begin{align*}
\big|\big|\Phi_t\big|\big|_{\Lambda^H_{2,\infty}} &=\sup_{\xi' \in {\Bbb R}^k} e^{-2 \pi^2 t{|\xi'|}^2} {\left( \int_{\R^{n-k}} e^{-4 \pi^2t{|\xi''|}^2} \, d\xi'' \right)}^{\frac{1}{2}} \\
&= {\left( \int_{\R^{n-k}} e^{-4 \pi^2 t{|\xi''|}^2} \, d\xi'' \right)}^{\frac{1}{2}} \\
&={\left( {(2 \pi \sqrt{t})}^{-\frac{n-k}{2}} \int_{\R^{n-k}} e^{-{|z''|}^2}\, dz'' \right)}^{\frac{1}{2}}={(4 \pi t)}^{-\frac{n-k}{4}}.
\end{align*}

\vskip.25in

\subsection{The wave equation: restriction to subspaces} 

\vskip.125in

We now consider  solutions to the Cauchy problem for the wave equation,
\begin{equation}\label{wave}
u_{tt}= \bigtriangleup u, \, (x,t)\in\R^{n}\times\R_+,\, \ u(x,0)=0, \ u_t(x,0)=F(x,y)\in L^2(\R^{n}).
\end{equation}

\noindent As is well-known, the map $F(x)\longrightarrow u(x,t)$ is bounded between a variety of function spaces; see, e.g., \cite{So93,St93} and the references  there. We show here  that for fixed $t\ne 0$ the solution $u(x,t)$ restricts in a natural way to any $k$-plane. 

By the standard integral representation of the solution to (\ref{wave}),
$$ u(x,t)=F*K_t(x):=c_n {\left( \frac{1}{t} \frac{d}{dt} \right)}^{\frac{n-3}{2}} \left\{ t^{n-2} A_t F(x,y) \right\},$$ 
with, as before,
$$ A_t F(x,y)=\int_{\mathbb S^{n-1}} F(x-tu)\, d\sigma(u),$$
where $d\sigma$ is the surface measure on the unit sphere in ${\Bbb R}^n$ and $K_t$ is the (fixed time) fundamental solution. The following is an immediate consequence of  Thm. \ref{bilinearnewwithattitude}.

\begin{theorem} \label{waveequationsolutionrestricted} Let $H$ be a $k$-dimensional affine subspace of ${\Bbb R}^{n}$. If $s<1-\frac{n-k}{2}$, then,
for all $t>0$,

$$ {\big|\big|{(I-\bigtriangleup_{x})}^{\frac{s}{2}} u(\cdot,t)|_{H}\big|\big|}_{L^2(H)} \leq C {\big|\big|F\big|\big|}_{L^2({\Bbb R}^{n})},$$ 
with the optimal constant  given by
$$  C_{s,n-k,t}={\big|\big|(I-\bigtriangleup_{x})^\frac{s}2K_t\big|\big|}_{\Lambda^{{H}}_{2,\infty}({\Bbb R}^{n})}.$$
\end{theorem}

\vskip.125in

\noindent{\bf Proof.} Since 
$\big| \big((I-\bigtriangleup_x)^sK_t\big)^{\widehat{\, \,}}(\xi)\big|\lesssim (1+|\xi|)^{s-1}$,  after putting $H$ into the same  form $\{(x',0)\in\R^n\, |\, x'\in\R^k\}$ as used above, we see that, for  $s<1$,
 $$\big|\big|(I-\bigtriangleup_x)^sK_t\big|\big|_{\Lambda^H_{2,\infty}}<\infty\quad \hbox{ iff }\quad  \int_{\R^{n-k}} (1+|\xi''|^2)^{2s-2} \, d\xi'' <\infty,$$
 which holds off $2s-2<-(n-k)$, i.e., $s<1-\frac{n-k}2$.

\subsection{The wave equation: product of solutions} 

Consider solutions $u$ and $v$ on $\Bbb R^{3+1}$ to
\smallskip
$$ u_{tt}=\bigtriangleup u, \ u(x,0)=0, \ u_t(x,0)=f(x),\quad\hbox{ and }\quad
v_{tt}=\bigtriangleup v, \ v(x,0)=0, \ v_t(x,0)=g(x). $$
We have the following estimate for products of solutions to these; see \cite{FK00} for similar results with an additional average in time.

\begin{theorem} \label{klainerman} With the notation above, for any fixed $t$, 
\begin{equation} \label{productsobolevwave} 
{||u(\cdot,t)v(\cdot, t)||}_{L^3({\Bbb R}^3)} \le C
{||u(\cdot,t)v(\cdot, t)||}_{L^2_{\frac{1}{2}}({\Bbb R}^3)} \leq C' {||f||}_{L^2({\Bbb R}^3)} \cdot {||g||}_{L^2({\Bbb R}^3)}, \end{equation} 
\begin{equation} \label{waveLpimproving} {||u(\cdot,t)v(\cdot,t)||}_{L^{p'}({\Bbb R}^3)} \leq C_{p,t} {||f||}_{L^p({\Bbb R}^3)} \cdot {||g||}_{L^p({\Bbb R}^3)} \ \text{for} \ \frac{5}{3} < p \leq 2. \end{equation} 
\end{theorem} 

\begin{remark} One then obtains $L^p\times L^p\to L^q$ estimates for $(\frac1p,\frac1q)$ in a nonclosed triangle with vertices at $(\frac12,\frac12), (\frac12,\frac13)$ and $(\frac35,\frac25)$.
\end{remark}

To prove Thm. \ref{klainerman}, observe that 

$$ u(x,t)=\frac{t}{4\pi}  \int_{S^2} f(x-ty) d\sigma(y),$$ 
so that estimate (\ref{waveLpimproving})   follows immediately from Cor. \ref{productsphere}. On the other hand,  (\ref{productsobolevwave}) follows from Thm. \ref{laceysobolevproduct}. 
In fact, since $ |\widehat{\sigma}(t \xi)| \leq C{(1+t |\xi|)}^{-1}$,  
 it  follows that 
\begin{equation} \label{waveu} {||u(\cdot,t)||}_{L^2_1({\Bbb R}^3)} \leq C{||f||}_{L^2({\Bbb R}^3)} \hbox{ and } 
 {||v(\cdot,t)||}_{L^2_1({\Bbb R}^3)} \leq C{||g||}_{L^2({\Bbb R}^3)},\end{equation} 
 with $C$ independent of $t$.
Let $r=s=1$, so that $\gamma=2-\frac{3}{2}=\frac{1}{2}$. Thm. \ref{laceysobolevproduct}, followed by (\ref{waveu}),  implies
$$ {||u(\cdot,t)v(\cdot,t)||}_{L^2_{\frac{1}{2}}({\Bbb R}^3)} \leq C {||u||}_{L^2_1({\Bbb R}^3)}  {||v||}_{L^2_1({\Bbb R}^3)} 
\leq C {||f||}_{L^2({\Bbb R}^3)} \cdot {||g||}_{L^2({\Bbb R}^3)}.$$ 
One then applies Sobolev embedding in $\R^3$ to obtain the first inequality in (\ref{productsobolevwave}).

\vskip.25in

\section{Bilinear oscillatory integral operators} 
\label{oscillatorysection}

\vskip.125in 

Let $\phi: {\Bbb R}^d \times {\Bbb R}^d \to {\Bbb R}$ be a smooth phase function and $\psi\in C_0^\infty(\Bbb R^{2d})$ a fixed amplitude. Define 

$$ T^{\lambda}_{\phi}F(x,y)=\int \int F(x-u,y-v) e^{2 \pi i \lambda \phi(u,v)} \psi(u,v) dudv. $$ 
If $\phi$ is \emph{nondegenerate}, i.e., the determinant of the Hessian matrix of $\phi$ does not vanish on the support of $\psi$, then it is a result of H\"ormander that 
$$ T^{\lambda}_{\phi}: L^2({\Bbb R}^{2d}) \to L^2({\Bbb R}^{2d}) \ \text{with norm} \  \le C \lambda^{-d}.$$ 
See \cite{St93} and the references  there for  this and related results. 
What we describe above works equally well in any ${\Bbb R}^n$, but this formulation leads  naturally to the context of trace inequalities, which we now present.  

\begin{definition} The \emph{bilinear surface associated with $\phi$ }  is

$$ S_{\phi}=\{(u,v,z, \delta_z \phi(u,v)): (u,v) \in U; z \in Z \}\subset \Bbb R^{3d+1},$$ 
\medskip

\noindent where $U:= supp(\psi)$, $Z=\pi_1(U)-\pi_1(U)$, with $\pi_1(u,v)=u$, and 

\begin{equation}\label{deltadef}
\delta_z \phi(u,v)=\phi(u,v)-\phi(u-z,v-z).
\end{equation}
\end{definition} 

\medskip

Let $\sigma_{\phi}$ denote the induced surface measure on $S_{\phi}$. 
We shall see below that  Thm. \ref{bilinearnewwithattitude} implies the following result. 

\begin{theorem} \label{generaloscillatory} Let $H$ be a $d$-dimensional plane in ${\Bbb R}^{2d}$. Then 

$$ {||(T^{\lambda}_{\phi}F)|_H||}_{L^2(H)} \leq   \sqrt{ \sup_{\xi \in {\Bbb R}^d} |\widehat{\sigma}_{\phi}(\xi, \lambda)|} \cdot {||F||}_{L^2({\Bbb R}^{2d})}.$$ 

\end{theorem} 

\vskip.125in 

We now apply our results to bilinear oscillatory integral operators, defined by

$$ M^{\lambda}_{\phi}(f,g)(x)=\int \int f(x-u) g(x-v) e^{2 \pi i \lambda \phi(u,v)} \psi(u,v) dudv.$$

\vskip.125in 

We have the following bound on the $L^2$-norm of $M^{\lambda}_{\phi}$, which follows at once from Thm. \ref{generaloscillatory}.  

\begin{corollary} \label{oscillatory} Let $M^{\lambda}_{\phi}$ be as above. Then 

\begin{equation} \label{oscillatoryest} {||M^{\lambda}_{\phi}(f,g)||}_2 \leq \sqrt{ \sup_{\xi \in {\Bbb R}^d} |\widehat{\sigma}_{\phi}(\xi, \lambda)|} \cdot {||f||}_2 {||g||}_2. \end{equation} 

\end{corollary} 

\vskip.125in

\begin{corollary} \label{oscillatoryhigherd} Suppose that the determinant of the Hessian matrix of $\phi$ does not vanish on the support of $\psi$. Then 

\begin{equation} \label{oscillatoryesthigherd} {||M^{\lambda}_{\phi}(f,g)||}_2 \leq C{\lambda}^{-\frac{d}{2}} {||f||}_2 {||g||}_2. \end{equation} 

\end{corollary}

\vskip.125in 

\begin{remark} \label{sharpness} The power of $\lambda$ in the estimate (\ref{oscillatoryesthigherd}) cannot, in general, be improved. This can be checked by taking $\phi(u,v)={|u|}^2+{|v|}^2$ and running a simple scaling argument. \end{remark} 

\vskip.125in 

\begin{example} \label{dotproduct} Let $\phi(u,v)=u \cdot v$. Then the conditions of Cor. \ref{oscillatoryhigherd} are satisfied and the conclusion holds. 
\end{example} 

\vskip.125in 

\noindent{\bf Proof of Theorem \ref{generaloscillatory}.} By Thm. \ref{bilinearnewwithattitude}, we need to compute the square root of 

$$ \sup_{\xi \in {\Bbb R}^{d}} \int {\left| \int \int e^{-2 \pi i \left( \frac{\xi-\eta}{2} \cdot u+ \frac{\xi+\eta}{2} \cdot v+\lambda \phi(u,v) \right)} \psi(u,v) \,dudv \right|}^2 d\eta$$

$$= \sup_{\xi \in {\Bbb R}^{d}} \int \int \dots \int e^{-2 \pi i \left( \frac{\xi-\eta}{2} \cdot (u-u')+ \frac{\xi+\eta}{2} \cdot (v-v')+\lambda (\phi(u,v)-\phi(u',v')) \right)} \psi(u,v) \psi(u',v') \,dudvdu'dv'  d\eta$$

$$=\sup_{\xi \in {\Bbb R}^{d}} \int \int \dots \int e^{-2 \pi i \left( \frac{\xi}{2} \cdot (u-u')+ \frac{\xi}{2} \cdot (v-v')+\lambda (\phi(u,v)-\phi(u',v')) \right)} \psi(u,v) \psi(u',v') \delta_0((u-u')-(v-v'))\, dudvdu'dv',$$ 
since the integral with respect to $\eta$ in the previous line yields $\delta_0((u-u')-(v-v'))$, where $\delta_0$ denotes the $\delta$-distribution at the origin. Making the change of variables $z=u-u', u=u, v=v$, we obtain

$$\sup_{\xi \in {\Bbb R}^{d}} \int \int \int e^{-2 \pi i (\xi \cdot z+\lambda \delta_z \phi(u,v))} \psi(u,v) \psi(u-z,v-z) dudvdz=\sup_{\xi \in {\Bbb R}^d} \widehat{\sigma}_{\phi}(\xi, \lambda),$$ 
and the proof is complete. 

\vskip.125in 

\noindent{\bf Proof of Corollary \ref{oscillatoryhigherd}.} By Theorem \ref{generaloscillatory}, we must estimate 

$$ \sup_{\xi \in {\Bbb R}^d} \int \int \int e^{-2 \pi i \left( z \cdot \xi+\lambda \delta_z \phi(u,v) \right)}\, \psi(u,v)\, \psi(u-z,v-z)\,dudvdz,$$ and by the method of stationary phase (see, e.g., \cite{St93}) our claim would follow if we could show that the rank of the $3d\times 3d$ Hessian matrix of $\Phi(u,v,z)=\delta_z \phi(u,v)$ has rank $\ge 2d$. Since the argument is local, we can use the Morse lemma and establish the result with 
$$\phi(u,v)=u_1^2 \pm u_2^2 \pm \dots \pm u_d^2 \pm v_1^2 \pm v_2^2 \pm \dots \pm v_d^2,$$ 
for which
$$ \Phi(u,v,z)=\sum_{i=1}^d \pm \left({(z_i-u_i)}^2-u_i^2 \right) +\sum_{j=1}^d \pm \left( {(z_j-v_j)}^2-v_j^2 \right), $$ 
and Cor. \ref{oscillatoryhigherd} follows. 

\vskip.25in

\section{Further thoughts and open problems} \label{thoughts}

\vskip.125in 

The key idea of this paper is that estimates for restrictions of convolutions (Theorem \ref{bilinearnewwithattitude}) provide a mechanism to obtain multilinear operator bounds, in the same way as the classical Young's inequality yields estimates for  linear operators. This approach does not rely on positivity, which was a drawback of the the technique in \cite{GGIP12}. 
Despite the applications we have presented in this paper, serious issues remain, especially regarding the bounds for bilinear operators arising in the geometric context. In \cite{GGIP12}, $L^2 \times L^2 \to L^1$ bounds were established for bilinear operators, with applications to Erd\H os/Falconer problems in geometry. These bounds relied heavily on the positivity of the kernel, which created significant obstacles, but the bounds were established in the arguably more natural range of exponents. While the problem of positivity was largely resolved in the current paper, the range of boundedness is centered around $L^2 \times L^2 \to L^2$. It would be extremely interesting to reconcile the methods in this paper and \cite{GGIP12} and come up with a unified set of bounds that do not rely on positivity of the kernel. We now state some concrete open problems. 

\begin{problem} It is proved in \cite{GGIP12} that if $B(f,g)$ is defined as above, then for $1 \leq r \leq 2$ and $K$ a non-negative finite measure, 
\begin{equation} \label{one} {||B(f,g)||}_{L^r({\Bbb R}^d)} \leq {||f||}_{L^2({\Bbb R}^d)} \cdot {||g||}_{L^2({\Bbb R}^d)} \cdot {\left( \int {\left| \widehat{K}(\xi, -\xi) \right|}^{r'} d\eta \right)}^{\frac{1}{r'}}. \end{equation}  
On the other hand, Cor. \ref{L2new} above shows that, if $2 \leq r \leq \infty$, 
\begin{equation} \label{two} {||B(f,g)||}_{L^r({\Bbb R}^d)} \leq {||f||}_{L^2({\Bbb R}^d)} \cdot {||g||}_{L^2({\Bbb R}^d)} 
\cdot {\left[ \int {\left( \int {\left|\widehat{K}\left(\frac{\xi-\eta}{2}, \frac{\xi+\eta}{2}\right)\right|}^2 d\eta \right)}^{\frac{r}{r-2}} d\xi \right]}^{\frac{r-2}{2r}}.\end{equation} 
The two estimates agree at $r=2$ and the positivity of $K$ is clearly irrelevant there. The question we ask is whether it is possible to reconcile (\ref{one}) and (\ref{two}) without assuming that $K$ is positive.
\end{problem} 

Our second problem attempts to further address the issue raised in Cor. \ref{restrictiontohyperplanes}. 

\begin{problem} It follows from (\ref{scramble}) that
$$ {\left( \int_{\R^d} {|\widehat{F}(x,x)|}^2 dx \right)}^{\frac{1}{2}} \leq {\left( \int_{\R^d} {\left[ \int_{\R^d} \left|F \left( \frac{x-y}{2}, \frac{x+y}{2} \right)
\right| dy \right]}^2 dx \right)}^{\frac{1}{2}}. $$

At least in two dimensions, it would be interesting to generalize this estimate to a universal $L^2$-restriction theorem where the left hand side is 
$$ {\left( \int {|\widehat{F}(x,\phi(x))|}^2 dx \right)}^{\frac{1}{2}},$$ with $\phi$ a suitably regular function of $x$, and the right hand side is a mixed norm depending on $\phi$. In the special case when $\phi(x)=x^2$, say, the right hand side should be comparable to the ${||F||}_{L^{\frac{6}{5}}({\Bbb R}^2)}$, consistent with the Stein-Tomas restriction theorem.

\end{problem}

%\newpage

\end{document}